\documentclass[a4paper,10pt]{article}
\usepackage{amsmath,amsthm,amssymb,amscd,epsfig,esint, mathrsfs,graphics,color}
 \usepackage{wrapfig}
\usepackage{verbatim}
\usepackage[english]{babel}
cm \oddsidemargin=.2true cm \topmargin=-.5truecm

\newcommand{\R}{\numberset{R}}

\theoremstyle{plain}
\newtheorem{thm}{Theorem}[section]

\newtheorem{proposition}[thm]{Proposition}
\newtheorem{lemma}[thm]{Lemma}

\theoremstyle{definition}
\newtheorem{definition}[thm]{Definition}

\def\Xint#1{\mathchoice
    {\XXint\displaystyle\textstyle{#1}}%
    {\XXint\textstyle\scriptstyle{#1}}%
    {\XXint\scriptstyle\scriptscriptstyle{#1}}%
    {\XXint\scriptscriptstyle\scriptscriptstyle{#1}}%
   \!\int}
\def\XXint#1#2#3{{\setbox0=\hbox{$#1{#2#3}{\int}$}
        \vcenter{\hbox{$#2#3$}}\kern-.5\wd0}}
\def\Mint{\Xint -}

\def\eps{\varepsilon}

\def\div{{\rm div}}

\def\R{\mathbb{R}}
\def\e{\varepsilon}
\def\dist{\mbox{dist }}
\numberwithin{equation}{section} \makeatletter

\renewcommand{\p@enumi}{\thesection.}

\title{\sc{Regularity results for solutions to a class of obstacle problems}
\footnotetext{\hspace{-0.35cm} 2010 \emph{Mathematics Subject
Classification}. 35J87, 49J40.
\endgraf
{\it Key words and phrases}.  }}
\author{  Michele Caselli, Andrea Gentile, Raffaella Giova}

\title{\textbf{Regularity results for solutions to obstacle problems with Sobolev coefficients}}
\begin{document}
\maketitle

\begin{abstract}
We establish the higher differentiability of solutions to a
class of obstacle problems of the type

\begin{equation*}
\min \left\{\int_{\Omega}f(x, Dv(x))dx\,:\, v\in
\mathcal{K}_\psi(\Omega)\right\},
\end{equation*}
where $\psi$ is a fixed function called obstacle,
$\mathcal{K}_\psi(\Omega)=\{v\in W^{1, p}_{\mathrm{loc}}(\Omega,
\R): v\ge\psi \text{ a.e. in }\Omega\}$ and the convex integrand $f$
satisfies $p$-growth conditions with respect to the gradient
variable. We derive that the higher differentiability property of
the weak solution $v$ is related to the regularity of the assigned
$\psi$, under a suitable Sobolev assumption on the partial map $x\mapsto D_\xi f(x, \xi)$. The main novelty is that such assumption is independent of the dimension $n$ and that,
in the case $p\le n-2$, improves previous known results.
\end{abstract}


\bigskip

\bigskip

\section{Introduction}

We are interested in the study of the regularity of the gradient of
the solutions to variational obstacle problems of the form

\begin{equation}\label{functionalobstacle}
\min \left\{\int_{\Omega}f(x, Dv(x))\,:\, v\in
\mathcal{K}_\psi(\Omega)\right\},
\end{equation}
where $\Omega\subset\R^n$ is a bounded open set, $\psi: \Omega
\mapsto [-\infty, +\infty)$ belonging to the Sobolev class $ W^{1,
\frac{p+2}{2}}_{\mathrm{loc}}$ is the \emph{obstacle}, and

$$\mathcal{K}_\psi(\Omega)=\{v\in W^{1, p}_{\mathrm{loc}}(\Omega, \R): v\ge\psi \text{ a.e. in }\Omega\}$$
is the class of the admissible functions.

\bigskip

Let us observe that $u\in W^{1,p}_{\mathrm{loc}}(\Omega)$ is a
solution to the obstacle problem \eqref{functionalobstacle} in
$\mathcal{K}_\psi(\Omega)$ if and only if $u \in
\mathcal{K}_\psi(\Omega)$ and $u$ is a solution to the variational
inequality

\begin{equation}\label{variationalinequality}
\int_{\Omega}\left<A(x, Du(x)), D(\varphi(x)-u(x))\right>dx\ge0\qquad\forall
\varphi\in \mathcal{K}_\psi(\Omega),
\end{equation}
where the operator $A(x, \xi): \Omega\times\R^n\to\R^n$ is defined
as follows

\begin{equation*}
A_i(x, \xi)=D_{\xi_i}f(x, \xi)\qquad\forall i=1,...,n.
\end{equation*}

We assume that $A$ is a $p$-harmonic operator, that satisfies the
following $p$-ellipticity and $p$-growth conditions with respect to
the $\xi$-variable. There exist positive constants $\nu, L, \ell$
and an exponent $p\ge2$ and a parameter $0\le \mu\le 1$  such that

$$\left<A(x, \xi)-A(x, \eta), \xi-\eta\right>\ge\nu|\xi-\eta|^2\left(\mu^2+|\xi|^2+|\eta|^2\right)^\frac{p-2}{2} \eqno{\rm (A1)}$$

%

$$\left|A(x, \xi)-A(x, \eta)\right|\le
L|\xi-\eta|\left(\mu^2+|\xi|^2+|\eta|^2\right)^\frac{p-2}{2}
\eqno{\rm (A2)}$$


$$
\left|A(x, \xi)\right|\le
\ell\left(\mu^2+|\xi|^2\right)^\frac{p-1}{2}, \eqno{\rm (A3)}$$
for all $\xi, \eta\in\R^n$ and for almost every $x\in\Omega.$\\

\bigskip

The regularity for solutions of obstacle problems has been object of
intense study not only in the case of variational inequalities
modelled upon the $p$- Laplacean energy \cite{Choe, CL, EP, MZ} but
also in the case of more general structures \cite{BLZ, BLOP, EHL, FM,Gavioli}\\
It is usually observed that the regularity of solutions to the
obstacle problems depends on the regularity of the obstacle itself:
for linear problems the solutions are as regular as the obstacle;
this is no longer the case in the nonlinear setting for general
integrands without any specific structure. Hence along the years, in this situation there has been an intense research activity in
which extra regularity has been imposed on the obstacle to balance
the nonlinearity (see \cite{BDM, BFM, CL, Fuchs, FM, Lindqvist})

\vspace{0.2cm} In some very recent papers the authors analyzed how
an extra differentiability of integer or fractional order of the
gradient of the obstacle transfers to the gradient of the solutions (see \cite{EP, EP2}).\\
The analysis comes from the fact that the regularity of the
solutions to the obstacle problem \eqref{functionalobstacle} is
strictly connected to the analysis of the regularity of the
solutions to partial differential equation of the form

\begin{equation}\label{star}
\div A(x,Du)= \div A(x, D\psi).
\end{equation}
It is well known that no extra differentiability properties for
the solutions can be expected even if the obstacle $\psi$ is smooth,
unless some assumption is given on the $x$-dependence of the operator $A$.
Therefore, inspired by recent results concerning the higher
differentiability of integer (\cite{EleMarMas, Gentile, Giova1, Giova2, GP, APdN1, APdN2}) and fractional (\cite{BCGOP, CGP}) order for the solutions
to elliptic equations or systems, in a number of papers the higher
differentiability of the solution of an obstacle problem is proved
under a suitable Sobolev assumption on the partial map $x \mapsto
A(x, \xi)$. More precisely, in \cite{EP} is proved the higher
differentiability of the solution of an homogeneous obstacle problem
with the energy density satisfying $p$-growth conditions; in
\cite{EP2} the integrand $f$ depends also on the $v$ variable; in
\cite{Gavioli} the energy density satisfies $(p,q)$-growth
conditions. The nonhomogeneous obstacle problem is considered
 in \cite{MaZ} when the energy density satisfies $p$-growth conditions and
 in \cite{CEP} when the energy density satisfies $(p,q)$-growth conditions.
 All previous quoted higher differentiability
results have been obtained under a $W^{1,r}$ with $r \le n$ Sobolev
assumption  on the dependence on $x$ of the operator $A$.\\
It is well known that the local boundedness of the solutions to a variational problem is a turning point in the regularity theory. Actually, in \cite{GP} it has been proved that, when dealing with bounded solutions to \eqref{star}, the higher differentiability holds true under weaker assumptions on the partial map $x\mapsto A(x, \xi)$ with respect to $W^{1, n}$. Recently, in \cite{CEP} it is proved that a local bound assumption on the obstacle $\psi$ implies a local bound for the solutions to the obstacle problem \eqref{functionalobstacle}, and this allows us to prove that the higher differentiability of solutions to \eqref{functionalobstacle} persists assuming that the partial map $x \mapsto A(x, \xi$ belongs to a Sobolev class that is not related to the dimension $n$ but to the growth exponent of the functional.\\
More precisely, we assume that there exists a non-negative function $\kappa\in
L^{p+2}_{\mathrm{loc}}(\Omega)$ such that

$$
\left|A(x, \xi)-A(y, \xi)\right|\le
\left(\kappa(x)+\kappa(y)\right)|x-y|\left(\mu^2+|\xi|^2\right)^\frac{p-1}{2}
\eqno{\rm (A4)}$$

\noindent for almost every $x, y\in\Omega$ and for all $\xi\in \R^n.$ The
condition (A4) is equivalent to assume that the operator $A$ has a
Sobolev-type dependence on the $x$-variable (see \cite{Hajlasz}).
Such assumption has been use for non constrained minimizers in \cite{KristensenMingione, KuusiMingione} 
\bigskip
We will prove a higher differentiability result assuming that $D
\psi\in W^{1, \frac{p+2}{2}}_{\mathrm{loc}}(\Omega).$ More
precisely, we shall prove the following.

\begin{thm}\label{thm1}
    Let $A(x, \xi)$ satisfy the conditions (A1)--(A4) for an exponent $p\ge2$ and
    let $u\in \mathcal{K}_\psi(\Omega)$ be a solution to the obstacle problem \eqref{variationalinequality}.
    Then, if $\psi \in L^{\infty}_\mathrm{loc}(\Omega)$ the
    following implication holds

$$D\psi\in W^{1, \frac{p+2}{2}}_{\mathrm{loc}}(\Omega) \,\, \Rightarrow   \,\, \left(\mu^2+\left|Du\right|^2\right)^\frac{p-2}{4}Du\in W^{1,2}_\mathrm{loc}(\Omega),$$

with the following estimate

\begin{align}\label{estimate}
\int_{B_{\frac{R}{4}}}&\left|D\left[\left(\mu^2+\left|Du\right|^2\right)^{\frac{p-2}{4}}Du\right]\right|^2dx
\le\notag\\
&\frac{c(\Arrowvert \psi\Arrowvert_{L^{\infty}}^2+\Arrowvert
	u\Arrowvert_{L^{p^*}\left(B_{R}\right)}^2)}{R^\frac{p+2}{2}}
\int_{B_{R}}\left[1+\left|D^2\psi\right|^\frac{p+2}{2}+
\left|D\psi\right|^\frac{p+2}{2}+\kappa^{p+2}+\left|Du\right|^p
\right] dx.
\end{align}

\end{thm}

Note that in the case $p<n-2$ Theorem \ref{thm1} improves the results in \cite{EP} and \cite{EP2}.
\noindent The proof of Theorem \ref{thm1} is achieved combining a
suitable a priori estimate for the second derivative of the local
solutions, obtained using the difference quotient method, with a
suitable approximation argument. The local boundedness allows us to
use an interpolation inequality that gives the higher local integrability
$L^{p+2}$ of the gradient of the solutions. Such higher
integrability is the key tool in order to weaken the assumption on
$\kappa$ that in previous results has been assumed at least in
$L^n$.\\
Moreover, our result is obtained under a weaker assumption also on
the gradient of the obstacle. Indeed, previous results assumed
$D\psi\in W^{1, p}$ while our assumption is $D\psi\in W^{1,
\frac{p+2}{2}}$ with $p>2$.\\
Finally, we observe that the assumption of boundedness of the
obstacle $\psi$ is needed to get the boundedness of the solution (see Theorem \ref{boundedness}). Therefore if we deal with a priori bounded
minimizers, then the result holds without the hypothesis $\psi \in
L^{\infty}$.

\section{Notations and preliminary results}

In this section we list the notations that we use in this paper and recall some tools that will be useful to prove our results.\\
We shall follow the usual convention and denote by $C$ or $c$ a general constant that may vary on different occasions, even within the same
line of estimates. Relevant dependencies on parameters and special
constants will be suitably emphasized using parentheses or
subscripts. All the norms we use on $\R^n$, $\R^N$ and $\R^{N\times n}$ will be
the standard Euclidean ones and denoted by $| \cdot |$ in all cases.
In particular, for matrices $\xi$, $\eta \in \R^{N\times n}$ we write $\langle
\xi, \eta \rangle : = \text{trace} (\xi^T \eta)$ for the usual inner
product of $\xi$ and $\eta$, and $| \xi | : = \langle \xi, \xi
\rangle^{\frac{1}{2}}$ for the corresponding Euclidean norm. When $a
\in \R^N$ and $b \in \R^n$ we write $a \otimes b \in \R^{N\times n}$ for the
tensor product defined
as the matrix that has the element $a_{r}b_{s}$ in its $r$-th row and $s$-th column.\\
For  a  $C^2$ function $f \colon \Omega\times\R^{N\times n} \to \R$, we write
$$
D_\xi f(x,\xi )[\eta ] := \frac{\rm d}{{\rm d}t}\Big|_{t=0} f(x,\xi
+t\eta )\quad \mbox{ and } \quad D_{\xi\xi}f(x,\xi )[\eta ,\eta ] :=
\frac{\rm d^2}{{\rm d}t^{2}}\Big|_{t=0} f(x,\xi +t\eta )
$$
for $\xi$, $\eta \in \R^{N\times n}$ and for almost every $x\in \Omega$.\\
With the symbol $B(x,r)=B_r(x)=\{y\in
\R^n:\,\, |y-x|<r\}$, we will denote the ball centered at $x$ of
radius $r$ and
$$(u)_{x_0,r}= \Mint_{B_r(x_0)}u(x)\,dx,$$
stands for the integral mean of $u$ over the ball $B_r(x_0)$. We
shall omit the dependence on the center  when it is clear from the context.
In the following, we will denote, for any ball
$B=B_r(x_0)=\{x\in\R^n: |x-x_0|<r\}\Subset\Omega$

\begin{equation}
\Mint_Bu(x)dx=\frac{1}{|B|}\int_Bu(x)dx.
\end{equation}

Here we recall some results that will be useful in the following.
The following Gagliardo-Niremberg type inequalities are stated as in
\cite{GP}. For the proofs see the Appendix A of \cite{CKP} and Lemma
3.5 in \cite{GiannettiPassa} (in case $p(x)\equiv p, \, \forall x$)
respectively.

\begin{lemma}\label{lemma5GP}
    For any $\phi\in C_0^1(\Omega)$ with $\phi\ge0$, and any $C^2$ map $v:\Omega\to\R^N$, we have

    \begin{align}\label{2.1GP}
    \int_\Omega&\phi^{\frac{m}{m+1}(p+2)}(x)|Dv(x)|^{\frac{m}{m+1}(p+2)}dx\notag\\
    \le&(p+2)^2\left(\int_\Omega\phi^{\frac{m}{m+1}(p+2)}(x)|v(x)|^{2m}dx\right)^\frac{1}{m+1}\cdot\left[\left(\int_\Omega\phi^{\frac{m}{m+1}(p+2)}(x)\left|D\phi(x)\right|^2\left|Dv(x)\right|^pdx\right)^\frac{m}{m+1}\right.\notag\\
    &\left.+n\left(\int_\Omega\phi^{\frac{m}{m+1}(p+2)}(x)\left|Dv(x)\right|^{p-2}\left|D^2v(x)\right|^2dx\right)^\frac{m}{m+1}\right],
    \end{align}

    for any $p\in(1, \infty)$ and $m>1$. Moreover, for any $\mu\in[0,1]$

    \begin{align}\label{2.2GP}
        \int_{\Omega}&\phi^2(x)\left(\mu^2+\left|Dv(x)\right|^2\right)^\frac{p}{2}\left|Dv(x)\right|^2dx\notag\\
        \le&c\Arrowvert v\Arrowvert_{L^\infty\left(\mathrm{supp}(\phi)\right)}^2\int_\Omega\phi^2(x)\left(\mu^2+\left|Dv(x)\right|^2\right)^\frac{p-2}{2}\left|D^2v(x)\right|^2dx\notag\\
        &+c\Arrowvert v\Arrowvert_{L^\infty\left(\mathrm{supp}(\phi)\right)}^2\int_\Omega\left(\phi^2(x)+\left|D\phi(x)\right|^2\right)\left(\mu^2+\left|Dv(x)\right|^2\right)^\frac{p}{2}dx,
    \end{align}

    for a constant $c=c(p).$
    \end{lemma}

By a density argument, one can easilt estimates \eqref{2.1GP} and \eqref{2.2GP} are still true for any map $v\in W^{2,p}_{\mathrm{loc}}(\Omega)$.\\

Moreover, if we recall Theorem 1.1 in \cite{CEP} in the case $p=q$
that suits with our ellipticity and growth assumptions:

\begin{thm}\label{boundedness}
	Let u in $K_{\psi}(\Omega)$ be a solution of
	\eqref{differenceinequality} under the assumptions (A1) and (A2). If
	the obstacle $\psi \in L^{\infty}_{\mathrm{loc}}(\Omega)$, then $u
	\in L^{\infty}_{\mathrm{loc}}(\Omega)$ and the following estimate
	\begin{equation}\label{est boundedness}
	\sup_{B_{R/2}} |u| \le \left[ \sup |\psi|+ \left(
	\int_{B_R}|u|^{p^*} dx\right)\right]^{\gamma}
	\end{equation}
	holds for every ball $B_R \Subset \Omega$, for $\gamma(n,p)>0$ and
	$c=c(\ell, \nu, p, n)$,
\end{thm}

We will use the auxiliary function $V_p:\R^n\to\R^n$, defined as 

\begin{equation}\label{Vp}
V_p(\xi):=\left(\mu^2+|\xi|^2\right)^\frac{p-2}{4}\xi,
\end{equation}

\noindent for which the following estimates hold (see \cite{GM})

\begin{lemma}\label{lemma6GP}
Let $1<p<\infty$. There is a constant $c=c(n, p)>0$ such that

\begin{equation}\label{lemma6GPestimate1}
    c^{-1}\left(\mu^2+|\xi|^2+|\eta|^2\right)^\frac{p-2}{2}\le\frac{\left|V
_p(\xi)-V_p(\eta)\right|^2}{|\xi-\eta|^2}\le c\left(\mu^2+|\xi|^2+|\eta|^2\right)^\frac{p-2}{2},
    \end{equation}

\noindent for any $\xi, \eta\in\R^n.$
Moreover, for a $C^2$ function $g$, there is a constant $C(p)$ such that

\begin{equation}\label{lemma6GPestimate2}
C^{-1}\left|D^2g\right|^2\left(\mu^2+\left|Dg\right|^2\right)^\frac{p-2}{2}\le\left|D\left(V(Dg)\right)\right|^2\le C\left|D^2g\right|^2\left(\mu^2+\left|Dg\right|^2\right)^\frac{p-2}{2}
\end{equation}.
\end{lemma}

The next lemma can be proved using an iteration technique, and will
be needed in the following. Its proof can be found for example in
\cite[Lemma 6.1]{Giusti}.

\begin{lemma}[Iteration Lemma]\label{iteration}
    Let $h: [\rho, R]\to \R$ be a nonnegative bounded function, $0<\theta<1$, $A, B\ge0$ and $\gamma>0$. Assume that

    $$
    h(r)\le\theta h(d)+\frac{A}{(d-r)^\gamma}+B
    $$

    for all $\rho\le r<d\le R_0<R.$ Then

    $$
    h(\rho)\le \frac{cA}{(R_0-\rho)^\gamma}+cB,
    $$

    where $c=c(\theta, \gamma)>0$.
\end{lemma}

\subsection{Difference quotient}
\medskip
\noindent \noindent In order to get the regularity of the solutions
of the problem \eqref{functionalobstacle}, we shall use the
difference quotient method. We recall here the definition
and basic results.
\begin{definition}
Given $h\in\mathbb{R}$, for every function
$F:\mathbb{R}^{n}\to\mathbb{R}$ the finite difference operator is
defined by
$$
\tau_{h}F(x)=F(x+h)-F(x).
$$
\end{definition}
\par
We recall some properties of the   finite difference operator that
will be needed in the sequel. We start with the description of some
elementary properties that can be found, for example, in
\cite{Giusti}.

\bigskip

\begin{proposition}\label{findiffpr}

Let $F$ and $G$ be two functions such that $F, G\in
W^{1,p}(\Omega)$, with $p\geq 1$, and let us consider the set
$$
\Omega_{|h|}:=\left\{x\in \Omega : dist(x,
\partial\Omega)>|h|\right\}.
$$
Then
\begin{itemize}
\item[$(d1)$] $\tau_{h}F\in W^{1,p}(\Omega_{|h|})$ and
$$
D_{i} (\tau_{h}F)=\tau_{h}(D_{i}F).
$$
\item[$(d2)$] If at least one of the functions $F$ or $G$ has support contained
in $\Omega_{|h|}$ then
$$
\int_{\Omega} F\, \tau_{h} G\, dx =\int_{\Omega} G\, \tau_{-h}F\,
dx.
$$
\item[$(d3)$] We have
$$
\tau_{h}(F G)(x)=F(x+h )\tau_{h}G(x)+G(x)\tau_{h}F(x).
$$
\end{itemize}
\end{proposition}

\noindent The next result about finite difference operator is a kind
of integral version of Lagrange Theorem.
\begin{lemma}\label{le1} If $0<\rho<R$, $|h|<\frac{R-\rho}{2}$, $1 < p <+\infty$,
 and $F, DF\in L^{p}(B_{R})$ then
$$
\int_{B_{\rho}} |\tau_{h} F(x)|^{p}\ dx\leq c(n,p)|h|^{p}
\int_{B_{R}} |D F(x)|^{p}\ dx .
$$
Moreover
$$
\int_{B_{\rho}} |F(x+h )|^{p}\ dx\leq  \int_{B_{R}} |F(x)|^{p}\ dx .
$$
\end{lemma}

We conclude this section recaling this result, that is proved in \cite{Giusti}.

\begin{lemma}\label{Giusti8.2}
	Let $f:\R^n\to\R^N$, $f\in L^p(B_R)$ with $1<p<+\infty$. Suppose that there exist $\rho\in(0, R)$ and $M>0$ such that
	
	$$
	\sum_{s=1}^{n}\int_{B_\rho}|\tau_{s, h}f(x)|^pdx\le M^p|h|^p
	$$
	
	for every $h<\frac{R-\rho}{s}$. Then $f\in W^{1,p}(B_R, \R^N)$. Moreover
	
	$$
	\Arrowvert Df \Arrowvert_{L^p(B_\rho)}\le M.
	$$		
\end{lemma}

\bigskip

\section{Proof of the Theorem \ref{thm1}}

The proof of the theorem will be divided in two steps: in the first
one, we will establish the a priori estimate, while in the second
one we will conclude through an approximation argument.

\begin{proof}

{\bf Step 1: The a priori estimate.} Suppose that $u$ is a local
solution to the obstacle problem in $\mathcal{K}_\psi(\Omega)$ such
that $$Du \in W^{1,2}_\mathrm{loc}(\Omega)\text{\qquad and\qquad}\left
(\mu^2+\left|Du\right|^2\right)^\frac{p-2}{4}Du\in
W^{1,2}_\mathrm{loc}(\Omega)$$. By estimate \eqref{est boundedness}
and Lemma \ref{lemma5GP}, we also have $|Du|\in
L^{p+2}_\mathrm{loc}(\Omega)$. Note that the a priori assumption
$|Du|\in L^{p+2}_\mathrm{loc}(\Omega)$ implies that the variational
inequality \eqref{variationalinequality}, by a simple density
argument, holds
true for every $\varphi \in W^{1, \frac{p+2}{2}}$. \\
\noindent In order to choose suitable test functions $\varphi$ in
\eqref{variationalinequality} that involve the different quotient of
the solution and at the same time belong to the class of the
admissible functions $\mathcal{K}_\psi(\Omega)$, we proceed as done
in \cite{EP}.

\bigskip

Let us fix a ball $B_{ R}\Subset \Omega$ and arbitrary radii $\frac{
R}{2}<r<s<t<\lambda r< R$, with $1<\lambda<2$. Let us consider a cut
off function $\eta\in C^\infty_0(B_t)$ such that $\eta\equiv 1$ on
$B_s$ and  $|\nabla \eta|\le \frac{c}{t-s}$. From now on, with no
loss of generality, we suppose $R<1$.\\

Let $v\in W^{1,
p}_0(\Omega)$ be such that

    \begin{equation}\label{cond}
    u-\psi+\tau v\ge 0\qquad\forall \tau\in[0, 1],
    \end{equation}

\noindent  and observe that $\varphi:=u+\tau v\in
\mathcal{K}_\psi(\Omega)$ for all $\tau \in [0, 1]$, since
$\varphi=u+\tau v\ge \psi $. For $|h|<\frac{R}{4}$, we consider

    \begin{equation}\label{v1}
    v_1(x)=\eta^2(x)\left[(u(x+h)-\psi(x+h))-(u(x)-\psi(x))\right],
    \end{equation}

\noindent  so we have  $v_1\in W^{1, \frac{p+2}{2}}_0(\Omega)$, and,
for any $\tau\in[0,1]$, $v_1$ satisfies \eqref{cond}. Indeed, for
a.e. $x \in \Omega$ and for any $\tau\in[0,1]$

    \begin{eqnarray*}
    u(x)-\psi(x)+\tau v_1(x) &=&
    u(x)-\psi(x)+\tau\eta^2(x)\left[(u-\psi)(x+h)-(u-\psi)(x)\right]\cr\cr
    &=&\tau\eta^2(x)(u-\psi)(x+h)+(1-\tau\eta^2(x))(u-\psi)(x)\ge 0,
    \end{eqnarray*}

 \noindent   since $u\in \mathcal{K}_\psi(\Omega)$ and $0\le\eta\le1$.\\
    So we can use $\varphi=u+\tau v_1$ as a test function in inequality \eqref{variationalinequality},  thus getting

    \begin{equation}\label{3.3}
    0\le\int_{\Omega}\left<A(x, Du(x)), D\left[\eta^2(x)\left[(u-\psi)(x+h)-(u-\psi)(x)\right]\right]\right>dx.
    \end{equation}

In a similar way, we define

    \begin{equation}\label{v2}
    v_2(x)=\eta^2(x-h)\left[(u-\psi)(x-h)-(u-\psi)(x)\right],
    \end{equation}

and we have $v_2\in W^{1, \frac{p+2}{2}}_0(\Omega)$, and
\eqref{cond} still is satisfied for any $\tau\in[0,1]$, since

    \begin{eqnarray*}
    u(x)-\psi(x)+\tau v_2(x)&=&
    u(x)-\psi(x)+\tau\eta^2(x-h)\left[(u-\psi)(x-h)-(u-\psi)(x)\right]\cr\cr
    &=&\tau\eta^2(x)(u-\psi)(x-h)+(1-\tau\eta^2(x-h))(u-\psi)(x)\ge
    0.
    \end{eqnarray*}

By using in \eqref{variationalinequality} as test function
$\varphi=u+\tau v_2$, we get

    \begin{equation}
    0\le\int_{\Omega}\left<A(x, Du(x)), D\left[\eta^2(x-h)\left[(u-\psi)(x-h)-(u-\psi)(x)\right]\right]\right>dx,
    \end{equation}

and by means of a change of variable, we obtain

    \begin{equation}\label{3.5}
    0\le\int_{\Omega}\left<A(x+h, Du(x+h)), D\left[\eta^2(x)\left[(u-\psi)(x)-(u-\psi)(x+h)\right]\right]\right>dx.
    \end{equation}

Now we can add \eqref{3.3} and \eqref{3.5}, thus getting

    \begin{align*}
    0\le&\int_{\Omega}\left<A(x, Du(x)), D\left[\eta^2(x)\left[(u-\psi)(x+h)-(u-\psi)(x)\right]\right]\right>dx\\
    &+\int_{\Omega}\left<A(x+h, Du(x+h)), D\left[\eta^2(x)\left[(u-\psi)(x)-(u-\psi)(x+h)\right]\right]\right>dx,
    \end{align*}

that is

\begin{equation*}
    0\le\int_{\Omega}\left<A(x, Du(x))-A(x+h, Du(x+h)), D\left[\eta^2(x)\left[(u-\psi)(x+h)-(u-\psi)(x)\right]\right]\right>dx,
    \end{equation*}


which implies

    \begin{align*}
    0\ge&\int_{\Omega}\left<A(x+h, Du(x+h))-A(x, Du(x)), \eta^2(x)D\left[(u-\psi)(x+h)-(u-\psi)(x)\right]\right>dx\notag\\
    &+\int_{\Omega}\left<A(x+h, Du(x+h))-A(x, Du(x)),
    2\eta(x)D\eta(x)\left[(u-\psi)(x+h)-(u-\psi)(x)\right]\right>dx.
    \end{align*}

Previous inequality can be rewritten as follows

\begin{align}\label{differenceinequality3}
    0\ge&\int_{\Omega}\left<A(x+h, Du(x+h))-A(x+h, Du(x)),\eta^2(x)(Du(x+h)-Du(x))\right>dx\notag\\
    &-\int_{\Omega}\left<A(x+h, Du(x+h))-A(x+h, Du(x)), \eta^2(x)(D\psi(x+h)-D\psi(x))\right>dx\notag\\
    &+\int_{\Omega}\left<A(x+h, Du(x+h))-A(x+h, Du(x)), 2\eta(x)D\eta(x)\tau_h\left(u(x)-\psi(x)\right)\right>dx\notag\\
    &+\int_{\Omega}\left<A(x+h, Du(x))-A(x, Du(x)),\eta^2(x)(Du(x+h)-Du(x))\right>dx\notag\\
    &-\int_{\Omega}\left<A(x+h, Du(x))-A(x, Du(x)), \eta^2(x)(D\psi(x+h)-D\psi(x))\right>dx\notag\\
    &+\int_{\Omega}\left<A(x+h, Du(x))-A(x, Du(x)), 2\eta(x)D\eta(x)\tau_h\left(u-\psi\right)\right>dx\notag\\
    =:& \,I+II+III+IV+V+VI,
\end{align}

\noindent so we have

\begin{equation}\label{differenceinequality}
    I\le |II|+|III|+|IV|+|V|+|VI|.
\end{equation}

\noindent By the ellipticity assumption (A1), we get

\begin{equation}\label{I}
    I\ge\nu\int_{\Omega}\eta^2(x)|\tau_hDu(x)|^2\left(\mu^2+\left|Du(x+h)\right|^2+\left|Du(x)\right|^2\right)^\frac{p-2}{2}dx.
\end{equation}

\noindent By virtue of assumption (A2), using Young's inequality with
exponents $\left(2, 2\right)$, and then H\"{o}lder's inequality with
exponents $\left(\frac{p+2}{4}, \frac{p+2}{p-2}\right)$, by the
properties of $\eta$, we infer

\begin{align}\label{II}
    |II|\le& \, L\int_{\Omega}\eta^2(x)|\tau_hDu(x)|\left(\mu^2+|Du(x+h)|^2+|Du(x)|^2\right)^\frac{p-2}{2}|\tau_hD\psi(x)|dx\notag\\
    \le&\varepsilon\int_{\Omega}\eta^2(x)|\tau_hDu(x)|^2\left(\mu^2+|Du(x+h)|^2+|Du(x)|^2\right)^\frac{p-2}{2}dx\notag\\
    &+c_\varepsilon(L)\int_{\Omega}\eta^2(x)|\tau_hD\psi(x)|^2\left(\mu^2+|Du(x+h)|^2+|Du(x)|^2\right)^\frac{p-2}{2}dx\notag\\
    \le&\varepsilon\int_{\Omega}\eta^2(x)|\tau_hDu(x)|^2\left(\mu^2+|Du(x+h)|^2+|Du(x)|^2\right)^\frac{p-2}{2}dx\notag\\
    &+c_\varepsilon(L)\left(\int_{B_{t}}|\tau_hD\psi(x)|^\frac{p+2}{2}dx\right)^\frac{4}{p+2}\cdot\left(\int_{B_{\lambda r}}\left(\mu^{p+2}+|Du(x+h)|^{p+2}\right)dx\right)^\frac{p-2}{p+2}\notag\\
    \le&\varepsilon\int_{\Omega}\eta^2(x)|\tau_hDu(x)|^2\left(\mu^2+|Du(x+h)|^2+|Du(x)|^2\right)^\frac{p-2}{2}dx\notag\\
    &+c_\varepsilon(L)|h|^2\left(\int_{B_{\lambda r}}\left|D^2\psi(x)\right|^\frac{p+2}{2}dx\right)^\frac{4}{p+2}\cdot\left(\int_{B_{\lambda r}}\left(\mu^{p+2}+|Du(x+h)|^{p+2}\right)dx\right)^\frac{p-2}{p+2},
\end{align}

\noindent where we used Lemma \ref{le1}.
Similarly, by Young's and H\"{o}lder's inequality, by virtue of the
properties of $\eta$, and Lemma \ref{le1}, we can estimate the term $|III|$ as follows

\begin{align}\label{III}
|III|\le& \, 2L\int_{\Omega}\eta |D \eta||\tau_hDu(x)|\left(\mu^2+|Du(x+h)|^2+|Du(x)|^2\right)^\frac{p-2}{2}|\tau_h\left(u-\psi\right)|dx\notag\\
\le&\varepsilon\int_{\Omega}\eta^2(x)|\tau_hDu(x)|^2\left(\mu^2+|Du(x+h)|^2+|Du(x)|^2\right)^\frac{p-2}{2}dx\notag\\
&+\frac{c_\varepsilon(L)}{(t-s)^2}\int_{B_t\setminus B_s}\left(\mu^2+|Du(x+h)|^2+|Du(x)|^2\right)^\frac{p-2}{2}|\tau_h\left(u-\psi\right)|^2dx\notag\\
\le&\varepsilon\int_{\Omega}\eta^2(x)|\tau_hDu(x)|^2\left(\mu^2+|Du(x+h)|^2+|Du(x)|^2\right)^\frac{p-2}{2}dx\notag\\
&+\frac{c_\varepsilon(L)}{(t-s)^2}\cdot |h|^2\left(\int_{B_{\lambda
r}}\left(\mu^{p+2}+|Du(x)|^{p+2}\right)dx\right)^\frac{p-2}{p+2}\cdot
\left(\int_{B_{\lambda
r}}|D(u-\psi)(x)|^\frac{p+2}{2}dx\right)^\frac{4}{p+2}.
\end{align}

In order to estimate the term $|IV|$, 
we use assumption (A4), Young's inequality with exponents $(2,
2)$ and the properties of $\eta$, thus getting

\begin{align*}
|IV|\le&|h|\int_{\Omega}\eta^2(x)\left(\kappa(x+h)+\kappa(x)\right)\left(\mu^2+|Du(x)|^2\right)^\frac{p-1}{2}|\tau_hDu(x)|dx\notag\\
\le&\varepsilon\int_{\Omega}\eta^2(x)\left|\tau_hDu(x)\right|^2\left(\mu^2+|Du(x+h)|^2+|Du(x)|^2\right)^\frac{p-2}{2}dx\notag\\
&+c_\varepsilon|h|^2\int_{B_t}\left(\kappa(x+h)+\kappa(x)\right)^2\left(\mu^2+|Du(x)|^2\right)^\frac{p}{2}dx,
\end{align*}
and using H\"{o}lder's inequality with exponents
$\left(\frac{p+2}{2}, \frac{p+2}{p}\right)$, and the properties of
$\eta$, we have

\begin{align}\label{IV}
|IV|\le&\varepsilon\int_{\Omega}\eta^2(x)\left|\tau_hDu(x)\right|^2\left(\mu^2+|Du(x+h)|^2+|Du(x)|^2\right)^\frac{p-2}{2}dx\notag\\
&+c_\varepsilon|h|^2\left(\int_{B_{\lambda r}}\kappa^{p+2}(x)dx\right)^\frac{2}{p+2}\cdot\left(\int_{B_t}\left(\mu^{p+2}+|Du(x)|^{p+2}\right)dx\right)^\frac{p}{p+2}.
\end{align}

In order to estimate the term $|V|$, we use the condition (A4)
again, than H\"{o}lder's inequality with exponents $\left(p+2,
\frac{p+2}{p-1}, \frac{p+2}{2}\right)$, the properties of $\eta$,
and the properties of difference quotients of Sobolev functions, so
we get

\begin{align}\label{V}
|V|\le&|h|\int_{\Omega}\eta^2(x)\left(\kappa(x+h)+\kappa(x)\right)\left(\mu^2+|Du(x)|^2\right)^\frac{p-1}{2}\left|\tau_hD\psi(x)\right|dx\notag\\
\le&|h|\left(\int_{B_t}\left(\kappa(x+h)+\kappa(x)\right)^{p+2}dx\right)^\frac{1}{p+2}\cdot\left(\int_{B_t}\left(\mu^{p+2}+|Du(x)|^{p+2}\right)dx\right)^\frac{p-1}{p+2}\notag\\
&\cdot\left(\int_{B_t}\left|\tau_hD\psi(x)\right|^\frac{p+2}{2}dx\right)^\frac{2}{p+2}\notag\\
\le&|h|^2\left(\int_{B_{\lambda r}}\kappa^{p+2}(x)dx\right)^\frac{1}{p+2}\cdot\left(\int_{B_t}\left(\mu^{p+2}+|Du(x)|^{p+2}\right)dx\right)^\frac{p-1}{p+2}\notag\\
&\cdot\left(\int_{B_{\lambda
r}}\left|D^2\psi(x)\right|^\frac{p+2}{2}dx\right)^\frac{2}{p+2},
\end{align}
where we used the assumption $D\psi \in W^{1, \frac{p+2}{2}}$ and
first estimate of Lemma \ref{le1}.

 For what concerns the term $|VI|$, using the condition
(A4), the properties of $\eta$, H\"{o}lder's inequality with
exponents $\left(p+2, \frac{p+2}{p-1}, \frac{p+2}{2}\right)$, and
the properties of difference quotients of Sobolev functions, we have

\begin{align}\label{VI}
|VI|\le&2|h|\int_{\Omega}\eta(x)\left|D\eta(x)\right|\left(\kappa(x+h)+\kappa(x)\right)\left(\mu^2+|Du(x)|^2\right)^\frac{p-1}{2}
\left|\tau_h\left(u-\psi\right)\right|dx\notag\\
\le&\frac{c|h|}{t-s}\left(\int_{B_t}\left(\kappa(x+h)+\kappa(x)\right)^{p+2}dx\right)^\frac{1}{p+2}\cdot\left(\int_{B_t}\left(\mu^{p+2}+|Du(x)|^{p+2}\right)dx\right)^\frac{p-1}{p+2}\notag\\
&\cdot\left(\int_{B_t}\left|\tau_h\left(u-\psi\right)\right|^\frac{p+2}{2}dx\right)^\frac{2}{p+2}\notag\\
\le&\frac{c|h|^2}{t-s}\left(\int_{B_{\lambda r}}\kappa(x)^{p+2}dx\right)^\frac{1}{p+2}\cdot\left(\int_{B_t}\left(\mu^{p+2}+|Du(x)|^{p+2}\right)dx\right)^\frac{p-1}{p+2}\notag\\
&\cdot\left(\int_{B_{\lambda
r}}\left|D\left(u-\psi\right)(x)\right|^\frac{p+2}{2}dx\right)^\frac{2}{p+2}.
\end{align}

Plugging \eqref{I}, \eqref{II}, \eqref{III}, \eqref{IV}, \eqref{V} and \eqref{VI} into \eqref{differenceinequality}, and choosing $\varepsilon=\frac{\nu}{6}$, and reabsorbing the terms with the same integral of the right-hand side of \eqref{I}, we get

\begin{align}\label{Plugestimates}
\nu\int_{\Omega}&\eta^2(x)|\tau_hDu(x)|^2\left(\mu^2+\left|Du(x+h)\right|^2+\left|Du(x)\right|^2\right)^\frac{p-2}{2}dx\notag\\
\le&c|h|^2\left(\int_{B_{\lambda r}}\left|D^2\psi(x)\right|^\frac{p+2}{2}dx\right)^\frac{4}{p+2}\cdot\left(\int_{B_{\lambda r}}\left(\mu^{p+2}+|Du(x)|^{p+2}\right)dx\right)^\frac{p-2}{p+2}\notag\\
&+\frac{c|h|^2}{(t-s)^2}\cdot \left(\int_{B_{\lambda r}}\left(\mu^{p+2}+|Du(x)|^{p+2}\right)dx\right)^\frac{p-2}{p+2}\cdot\left(\int_{B_{\lambda r}}|D\left(u-\psi\right)(x)|^\frac{p+2}{2}dx\right)^\frac{4}{p+2}\notag\\
&+c|h|^2\left(\int_{B_{\lambda r}}\kappa^{p+2}(x)dx\right)^\frac{2}{p+2}\cdot\left(\int_{B_t}\left(\mu^{p+2}+|Du(x)|^{p+2}\right)dx\right)^\frac{p}{p+2}\notag\\
&+c|h|^2\left(\int_{B_{\lambda r}}\kappa^{p+2}(x)dx\right)^\frac{1}{p+2}\cdot\left(\int_{B_t}\left(\mu^{p+2}+|Du(x)|^{p+2}\right)dx\right)^\frac{p-1}{p+2}\notag\\
&\cdot\left(\int_{B_{\lambda r}}\left|D^2\psi(x)\right|^\frac{p+2}{2}dx\right)^\frac{2}{p+2}\notag\\
&+\frac{c|h|^2}{t-s}\left(\int_{B_{\lambda r}}\kappa(x)^{p+2}dx\right)^\frac{1}{p+2}\cdot\left(\int_{B_t}\left(\mu^{p+2}+|Du(x)|^{p+2}\right)dx\right)^\frac{p-1}{p+2}\notag\\
&\cdot\left(\int_{B_{\lambda
r}}\left|D\left(u-\psi\right)(x)\right|^\frac{p+2}{2}dx\right)^\frac{2}{p+2}.
\end{align}

Now we apply Young's inequality with exponents $\left(\frac{p+2}{4},
\frac{p+2}{p-2}\right)$ to the first two terms of the right-hand side of \eqref{Plugestimates}, Young's inequality with exponents $\left(\frac{p+2}{2},
\frac{p+2}{p}\right)$  to the third one, and $\left(p+2, \frac{p+2}{p-1},
\frac{p+2}{2}\right)$ to the last to terms, and since $u\in \mathcal{K}_{\psi}(\Omega)$,
we have

\begin{align}
\nu\int_{\Omega}&\eta^2(x)|\tau_hDu(x)|^2\left(\mu^2+\left|Du(x+h)\right|^2+\left|Du(x)\right|^2\right)^\frac{p-2}{2}dx\notag\\
\le&\varepsilon|h|^2\int_{B_{\lambda r}}\left(\mu^{p+2}+|Du(x+h)|^{p+2}\right)dx+c_\varepsilon|h|^2\int_{B_{\lambda r}}\left|D^2\psi(x)\right|^\frac{p+2}{2}dx\notag\\
&+\varepsilon|h|^2\int_{B_{\lambda r}}\left(\mu^{p+2}+|Du(x)|^{p+2}\right)dx+\frac{c_\varepsilon|h|^2}{(t-s)^{\frac{p+2}{2}}}\cdot\int_{B_{\lambda r}}\left|D\psi(x)\right|^\frac{p+2}{2}dx\notag\\
&+\varepsilon|h|^2\int_{B_{\lambda r}}\left(\mu^{p+2}+|Du(x)|^{p+2}\right)dx+c_\varepsilon|h|^2\int_{B_{\lambda r}}\kappa^{p+2}(x)dx\notag\\
&+\varepsilon|h|^2\int_{B_{\lambda r}}\left(\mu^{p+2}+|Du(x)|^{p+2}\right)dx
+c_\varepsilon|h|^2\int_{B_{\lambda r}}\left|D\psi(x)\right|^\frac{p+2}{2}dx\notag\\
&+\frac{c_\varepsilon|h|^2}{(t-s)^\frac{p+2}{2}}\int_{B_{\lambda r}}\kappa^{p+2}(x)dx.
\end{align}

Recalling the right-hand side of the inequality \eqref{lemma6GPestimate1} in Lemma \ref{lemma6GP}, we get

\begin{align}
\nu\int_\Omega&\eta^2(x)\left|\tau_hV_p\left(Du(x)\right)\right|^2dx\notag\\
\le&\varepsilon|h|^2\int_{B_{\lambda r}}\left(\mu^{p+2}+|Du(x+h)|^{p+2}\right)dx+c_\varepsilon|h|^2\int_{B_{\lambda r}}\left|D^2\psi(x)\right|^\frac{p+2}{2}dx\notag\\
&+\varepsilon|h|^2\int_{B_{\lambda r}}\left(\mu^{p+2}+|Du(x)|^{p+2}\right)dx+\frac{c_\varepsilon|h|^2}{(t-s)^{\frac{p+2}{2}}}\cdot\int_{B_{\lambda r}}\left|D\psi(x)\right|^\frac{p+2}{2}dx\notag\\
&+\varepsilon|h|^2\int_{B_{\lambda r}}\left(\mu^{p+2}+|Du(x)|^{p+2}\right)dx+c_\varepsilon|h|^2\int_{B_{\lambda r}}\kappa^{p+2}(x)dx\notag\\
&+\varepsilon|h|^2\int_{B_{\lambda r}}\left(\mu^{p+2}+|Du(x)|^{p+2}\right)dx
+c_\varepsilon|h|^2\int_{B_{\lambda r}}\left|D\psi(x)\right|^\frac{p+2}{2}dx\notag\\
&+\frac{c_\varepsilon|h|^2}{(t-s)^\frac{p+2}{2}}\int_{B_{\lambda r}}\kappa^{p+2}(x)dx.
\end{align}

Now we divide both sides by $|h|^2$ and use the Lemma \ref{Giusti8.2}, thus getting

\begin{align}\label{stima derivate seconde}
    \int_\Omega&\eta^2(x)\left|DV_p(Du(x))\right|^2dx\le4\varepsilon\int_{B_{\lambda r}}\left(\mu^{p+2}+|Du(x)|^{p+2}\right)dx+c_\varepsilon\int_{B_{\lambda r}}\left|D^2\psi(x)\right|^\frac{p+2}{2}dx\notag\\
    &+c_\varepsilon\int_{B_{\lambda r}}\left|D\psi(x)\right|^\frac{p+2}{2}dx+\frac{c_\varepsilon}{(t-s)^\frac{p+2}{2}}\int_{B_{\lambda r}}\left|D\psi(x)\right|^\frac{p+2}{2}dx+c_\varepsilon\int_{B_{\lambda r}}\kappa^{p+2}(x)dx\notag\\
    &+\frac{c_\varepsilon}{(t-s)^\frac{p+2}{2}}\int_{B_{\lambda r}}\kappa^{p+2}(x)dx
\end{align}

and, by left-hand side of inequality \eqref{lemma6GPestimate2},

\begin{align}\label{beforeInterpolation}
\int_{\Omega}&\eta^2(x)\left(\mu^2+\left|Du(x)\right|^2\right)^\frac{p-2}{2}\left|D^2u(x)\right|^2dx\le\int_\Omega\eta^2(x)\left|DV_p(Du(x))\right|^2dx\notag\\
\le&4\varepsilon\int_{B_{\lambda r}}\left(\mu^{p+2}+|Du(x)|^{p+2}\right)dx+c_\varepsilon\int_{B_{\lambda r}}\left|D^2\psi(x)\right|^\frac{p+2}{2}dx+c_\varepsilon\int_{B_{\lambda r}}\left|D\psi(x)\right|^\frac{p+2}{2}dx\notag\\
&+\frac{c_\varepsilon}{(t-s)^\frac{p+2}{2}}\int_{B_{\lambda r}}\left|D\psi(x)\right|^\frac{p+2}{2}dx+c_\varepsilon\int_{B_{\lambda r}}\kappa^{p+2}(x)dx+\frac{c_\varepsilon}{(t-s)^\frac{p+2}{2}}\int_{B_{\lambda r}}\kappa^{p+2}(x)dx.
\end{align}

By inequality \eqref{2.2GP} we have

\begin{align}
\int_{\Omega}&\eta^2(x)\left(\mu^2+\left|Du(x)\right|^2\right)^\frac{p}{2}\left|Du(x)\right|^2dx\notag\\
\le&c\Arrowvert u\Arrowvert_{L^\infty\left(\mathrm{supp}(\eta)\right)}^2\int_\Omega\eta^2(x)\left(\mu^2+\left|Du(x)\right|^2\right)^\frac{p-2}{2}\left|D^2u(x)\right|^2dx\notag\\
&+c\Arrowvert
u\Arrowvert_{L^\infty\left(\mathrm{supp}(\eta)\right)}^2\int_{\Omega}\left(|\eta(x)|^2+\left|D\eta(x)\right|^2\right)
\left(\mu^2+\left|Du(x)\right|^2\right)^\frac{p}{2}dx.
\end{align}

Hence, thanks to estimate \eqref{beforeInterpolation}, and the
properties of $\eta$ we infer

\begin{align}\label{beforeIteration1}
\int_{\Omega}&\eta^2(x)\left(\mu^2+\left|Du(x)\right|^2\right)^\frac{p}{2}\left|Du(x)\right|^2dx\le \varepsilon\cdot c\Arrowvert u\Arrowvert_{L^\infty\left(B_{\lambda r}\right)}^2\int_{B_{\lambda r}}\left(\mu^{p+2}+\left|Du(x)\right|^{p+2}\right)dx\notag\\
&+c_\varepsilon\Arrowvert u\Arrowvert_{L^\infty\left(B_{\lambda r}\right)}^2\int_{B_{\lambda r}}\left|D^2\psi(x)\right|^\frac{p+2}{2}dx+c_\varepsilon\Arrowvert u\Arrowvert_{L^\infty\left(B_{\lambda r}\right)}^2\int_{B_{\lambda r}}\left|D\psi(x)\right|^\frac{p+2}{2}dx\notag\\
&+\frac{c_\varepsilon\Arrowvert u\Arrowvert_{L^\infty\left(B_{\lambda r}\right)}^2}{(t-s)^\frac{p+2}{2}}\int_{B_{\lambda r}}\left|D\psi(x)\right|^\frac{p+2}{2}dx+c_\varepsilon\Arrowvert u\Arrowvert_{L^\infty\left(B_{\lambda r}\right)}^2\int_{B_{\lambda r}}\kappa^{p+2}(x)dx\notag\\
&+\frac{c_\varepsilon\Arrowvert u\Arrowvert_{L^\infty\left(B_{\lambda r}\right)}^2}{(t-s)^\frac{p+2}{2}}\int_{B_{\lambda r}}\kappa^{p+2}(x)dx+\frac{c_\varepsilon\Arrowvert u\Arrowvert_{L^\infty\left(B_{\lambda r}\right)}^2}{(t-s)^2}\int_{B_{\lambda r}}\left(\mu^2+\left|Du(x)\right|^2\right)^\frac{p}{2}dx.
\end{align}

Taking into account the properties of $\eta$ again, since $p\ge2$
and  $t-s<1$, we obtain

\begin{align*}
\int_{B_{r}}&\left(\mu^2+\left|Du(x)\right|^2\right)^\frac{p}{2}\left|Du(x)\right|^2dx\le
\varepsilon\cdot c\Arrowvert u\Arrowvert_{L^\infty\left(B_{R}\right)}^2\int_{B_{\lambda r}}\left(\mu^{p+2}+\left|Du(x)\right|^{p+2}\right)dx\notag\\
&+c_\varepsilon\Arrowvert u\Arrowvert_{L^\infty\left(B_{R}\right)}^2\int_{B_{R}}\left|D^2\psi(x)\right|^\frac{p+2}{2}dx
+c_\varepsilon\Arrowvert u\Arrowvert_{L^\infty\left(B_{R}\right)}^2\int_{B_{R}}\left|D\psi(x)\right|^\frac{p+2}{2}dx\notag\\
&+c_\varepsilon\Arrowvert
u\Arrowvert_{L^\infty\left(B_{R}\right)}^2\int_{B_{R}}\kappa^{p+2}(x)dx
+\frac{c_\varepsilon\Arrowvert
u\Arrowvert_{L^\infty\left(B_{R}\right)}^2}{(t-s)^\frac{p+2}{2}}
\left[\int_{B_{R}}\left|D\psi(x)\right|^\frac{p+2}{2}dx\right.\notag\\
&\left.+\int_{B_{R}}\kappa^{p+2}(x)dx+\int_{B_{R}}\left(\mu^2+\left|Du(x)\right|^2\right)^\frac{p}{2}dx\right],
\end{align*}

and choosing $\varepsilon$ such that $\varepsilon\cdot c\Arrowvert
u\Arrowvert_{L^\infty\left(B_{R}\right)}^2\le\frac{1}{2}$, previous
estimate becomes

\begin{align}\label{beforeIteration3}
\int_{B_{r}}&\left|Du(x)\right|^{p+2}dx\le\int_{B_{r}}\left(\mu^2+\left|Du(x)\right|^2\right)^\frac{p}{2}\left|Du(x)\right|^2dx\le
\frac{1}{2}\int_{B_{\lambda r}}\left|Du(x)\right|^{p+2}dx\notag\\
&+c\Arrowvert
u\Arrowvert_{L^\infty\left(B_{R}\right)}^2\left[\int_{B_{R}}\left|D^2\psi(x)\right|^\frac{p+2}{2}dx+
\int_{B_{R}}\left|D\psi(x)\right|^\frac{p+2}{2}dx +\int_{B_{R}}\kappa^{p+2}(x)dx\right]\notag\\
&+\frac{c\Arrowvert
u\Arrowvert_{L^\infty\left(B_{R}\right)}^2}{(t-s)^\frac{p+2}{2}}\left[\int_{B_{R}}\left|D\psi(x)\right|^\frac{p+2}{2}dx
+\int_{B_{R}}\kappa^{p+2}(x)dx+\int_{B_{R}}\left|Du(x)\right|^p dx
+c( \mu,p)|B_{R}|\right]+c( \mu,p)|B_{R}|,
\end{align}

%

\noindent where $c=c(p, L, \nu, \mu)$ is independent of $t$ and $s$. Since
\eqref{beforeIteration3} is valid for any $\frac{R}{2}<r<s<t<\lambda
r<R<1$, taking the limit as $s\to r$ and $t\to\lambda r$, we get

\begin{align}\label{beforeIteration5}
\int_{B_{r}}&\left|Du(x)\right|^{p+2}dx\le \frac{1}{2}\int_{B_{\lambda r}}\left|Du(x)\right|^{p+2}dx\notag\\
&+c\Arrowvert
u\Arrowvert_{L^\infty\left(B_{R}\right)}^2\left[\int_{B_{R}}\left|D^2\psi(x)\right|^\frac{p+2}{2}dx+\int_{B_{R}}\left|D\psi(x)\right|^\frac{p+2}{2}dx
+\int_{B_{R}}\kappa^{p+2}(x)dx\right]\notag\\
&+\frac{c\Arrowvert
u\Arrowvert_{L^\infty\left(B_{R}\right)}^2}{r^\frac{p+2}{2}(\lambda-1)^\frac{p+2}{2}}\left[\int_{B_{R}}\left|D\psi(x)\right|^\frac{p+2}{2}dx
+\int_{B_{R}}\kappa^{p+2}(x)dx+\int_{B_{R}}\left|Du(x)\right|^p dx
+c( \mu,p)|B_{R}|\right]+c( \mu,p)|B_{R}|.
\end{align}

 Now, setting

\begin{equation*}
 h(r)=\int_{B_{r}}\left|Du(x)\right|^{p+2}dx,
\end{equation*}

\begin{equation*}
 A=c\Arrowvert u\Arrowvert_{L^\infty\left(B_{R}\right)}^2\left[\int_{B_{R}}\left|D\psi(x)\right|^\frac{p+2}{2}dx+
 \int_{B_{R}}\kappa^{p+2}(x)dx+\int_{B_{R}}\left|Du(x)\right|^p dx +c(
\mu,p)|B_{R}|\right],
\end{equation*}

and

\begin{equation*}
B=c\Arrowvert
u\Arrowvert_{L^\infty\left(B_{R}\right)}^2\left[\int_{B_{R}}\left|D^2\psi(x)\right|^\frac{p+2}{2}dx+
\int_{B_{R}}\left|D\psi(x)\right|^\frac{p+2}{2}dx+\int_{B_{R}}\kappa^{p+2}(x)dx\right]+c(
\mu,p)|B_{R}|,
\end{equation*}

we can use Lemma \ref{iteration}, with

\begin{equation*}
\theta=\frac{1}{2}\qquad\text{ and }\qquad\gamma=\frac{p+2}{2},
\end{equation*}

thus obtaining

\begin{align}\label{Iteration}
\int_{B_{\frac{R}{2}}}&\left|Du\right|^{p+2} dx\le c\Arrowvert
u\Arrowvert_{L^\infty\left(B_{R}\right)}^2\left[\int_{B_{R}}\left|D^2\psi\right|^\frac{p+2}{2}dx+\int_{B_{R}}\left|D\psi\right|^\frac{p+2}{2}dx\right.\notag\\
&\left.+\int_{B_{R}}\kappa^{p+2}dx\right]+c(
\mu,p)|B_{R}|+\frac{c\Arrowvert
u\Arrowvert_{L^\infty\left(B_{R}\right)}^2}{R^\frac{p+2}{2}}\left[\int_{B_{R}}\left|D\psi\right|^\frac{p+2}{2}dx+\int_{B_{R}}\kappa^{p+2}dx\right.\notag\\
&\left.+\int_{B_{R}}\left|Du\right|^p dx+c( \mu,p)|B_{R}|\right].
\end{align}
Since $R<1$, estimate \eqref{Iteration} can be written as follows

\begin{equation}\label{est pulita}
\int_{B_{\frac{R}{2}}}\left|Du\right|^{p+2} dx\le \frac{c\Arrowvert
u\Arrowvert_{L^\infty\left(B_{R}\right)}^2}{R^\frac{p+2}{2}}\int_{B_{R}}\left[1+\left|D^2\psi\right|^\frac{p+2}{2}+
\left|D\psi\right|^\frac{p+2}{2}+\kappa^{p+2}+\left|Du\right|^p
\right] dx.
\end{equation}

\noindent Now, we consider the estimate in \eqref{stima derivate
seconde} choosing a cut off function $\eta\in
C^\infty_0(B_{\frac{R}{2}})$ such that $\eta\equiv 1$ on
$B_{\frac{R}{4}}$; so that, thanks to \eqref{est pulita}, we
obtain

\begin{align*}
\int_{B_{\frac{R}{4}}}&\left|DV_p(Du(x))\right|^2 dx\le
\frac{c\Arrowvert
u\Arrowvert_{L^\infty\left(B_{R}\right)}^2}{R^\frac{p+2}{2}}\int_{B_{R}}\left[1+\left|D^2\psi\right|^\frac{p+2}{2}+
\left|D\psi\right|^\frac{p+2}{2}+\kappa^{p+2}+\left|Du\right|^p
\right] dx.
\end{align*}

By virtue of estimate \eqref{est boundedness}, we conclude
with

\begin{align}\label{final estimate }
 \int_{B_{\frac{R}{4}}}\left|DV_p(Du(x))\right|^2 dx
\le \frac{c(\Arrowvert \psi\Arrowvert_{L^{\infty}}^2+\Arrowvert
u\Arrowvert_{L^{p^*}\left(B_{R}\right)}^2)}{R^\frac{p+2}{2}}
\int_{B_{R}}\left[1+\left|D^2\psi\right|^\frac{p+2}{2}+
\left|D\psi\right|^\frac{p+2}{2}+\kappa^{p+2}+\left|Du\right|^p
\right] dx.
\end{align}
%

\bigskip

{\bf Step 2: The approximation.}

Fix a compact set  $\Omega'\Subset \Omega$,  and for a smooth kernel
$\phi \in C^{\infty}_{c}(B_{1}(0))$ with $\phi \geq 0$ and
$\int_{B_{1}(0)} \! \phi = 1$, let us consider the corresponding
family of mollifiers $( \phi_{\eps})_{\eps >0}$ and put

$$\kappa_\e=\kappa\ast\phi_\e, \qquad \qquad \psi_\e=\psi\ast\phi_\e,$$

$$ \mathcal{K}_{\psi_\e}(\Omega)=\{v\in u+W^{1,
p}_0(\Omega): v\ge\psi_\e \text{ almost everywhere in }\Omega\}$$
and
 \begin{equation}\label{A_epsilon}
A_\varepsilon(x,\xi)=\int_{B_1}\phi(\omega)A(x+\varepsilon\omega,\xi)\,\mathrm{d}\omega
\end{equation}
on $\Omega'$, for each positive $\eps < \dist (\Omega',\Omega)$. The
assumptions (A1)--(A3) imply that
$$\langle A_\varepsilon(x,\xi)-A_\varepsilon(x,\eta),\xi-\eta\rangle \ge    \nu|\eta-\xi|^2(\mu^2+|\xi|^2+|\eta|^2)^{\frac{p-2}{2}}
\eqno{\rm (A1^\prime)}$$
$$|A_\varepsilon(x,\xi)-A_\varepsilon(x,\eta)|\le L |\xi-\eta |(\mu^2+|\xi|^2+|\eta|^2)^{\frac{p-2}{2}}\eqno{\rm (A2^\prime)}$$
$$
\left|A_\e(x, \xi)\right|\le
\ell\left(\mu^2+|\xi|^2\right)^\frac{p-1}{2}, \eqno{\rm
(A3^\prime)}$$ By virtue of assumption (A4), we have that
$$| A_\varepsilon(x,\xi)-A_\varepsilon(y,\xi)|\le  (\kappa_{\e}(x)+\kappa_\e(y))|x-y|(1+| \xi |^2)^{\frac{p-1}{2}}\,.\eqno{\rm (A4^\prime)}$$
for almost every $x,y\in\Omega$ and for all $\xi,\eta \in
\mathbb{R}^{n}$. Let $u$ be a solution of the variational inequality
\eqref{variationalinequality} and let fix a ball $B_R\Subset
\Omega'$. Let us denote by $u_\eps \in W^{1,p}(B_R)$ the solution of
the problem

\begin{equation}\label{variationalinequality2}
\int_{\Omega}\left<A_\e(x, Dv),
D(\varphi-v)\right>dx\ge0\qquad\forall \varphi\in
\mathcal{K}_{\psi_\e}(\Omega).
\end{equation}
Thanks to \cite[Theorem 1.1]{EP} we have
$\left(\mu^2+\left|Du_\e\right|^2\right)^\frac{p-2}{4}Du_\e\in
W^{1,2}_\mathrm{loc}(B_R) $ and, since $A_\eps$ satisfies conditions
(A1$^\prime$)--(A4$^\prime$), for $\eps$ sufficiently small, we are
legitimate to apply estimate \eqref{final estimate } to get

\begin{align}\label{final estimate_eps }
\int_{B_{\frac{r}{4}}}\left|DV_p(Du_\e(x))\right|^2 dx\le
\frac{c(\Arrowvert \psi_\e\Arrowvert_{L^{\infty}}^2+\Arrowvert
u_\e\Arrowvert_{L^{p^*}\left(B_{r}\right)}^2)}{r^\frac{p+2}{2}} &
\int_{B_{r}}\left[1+\left|D^2\psi_\e\right|^\frac{p+2}{2}+
\left|D\psi_\e\right|^\frac{p+2}{2}+\kappa_\e^{p+2}+\left|Du_\e\right|^p
\right] dx.
\end{align}
for every ball $B_r\Subset B_R$ and for a constant $c=c()$.

\noindent We recall that, since $D\psi \in W^{1,
\frac{p+2}{2}}_\mathrm{loc}(\Omega)$ and $\kappa\in
L^{p+2}_\mathrm{loc}(\Omega)$, then

\begin{equation}\label{convpsi}
   D\psi_\eps \to D\psi \qquad  \mathrm{and}\qquad D^2\psi_\eps \to D^2\psi \qquad \mathrm{strongly\,\,in}\,\,
   L^{\frac{p+2}{2}}_{\mathrm{loc}}(\Omega'),
\end{equation}

\begin{equation}\label{convk}
   \kappa_\eps \to \kappa\qquad \mathrm{strongly\,\,in}\,\, L^{p+2}_{\mathrm{loc}}(\Omega').
\end{equation}

Since from (A3$^\prime$) the function $|A_\varepsilon(x,Du)|\le \ell
(\mu^2+|Du|)^{p-1}$ and  since $A_\varepsilon(x,Du)$ converges
almost everywhere to $A(x,Du)$, by the dominated convergence Theorem
we have
\begin{equation}\label{convdf}
A_\varepsilon(x,Du)\to A(x,Du)\qquad
\mathrm{strongly\,\,in}\,\,L^{\frac{p}{p-1}}_{\mathrm{loc}}(\Omega').
\end{equation}

Using the ellipticity condition (A1$^\prime$) we have
\begin{eqnarray}\label{convforte}
&&\int_{B_R}(\mu^2+|Du|^2+|Du_\eps|^2)^{\frac{p-2}{2}}|Du_\eps-Du|^2\,dx\cr\cr
&\le&  \int_{B_R}\Big\langle A_\eps(x,Du_\eps)-A_\e(x,Du),
Du_\eps-Du\Big\rangle\,dx\cr\cr &=&\int_{B_R}\Big\langle
A_\eps(x,Du_\eps), Du_\eps-Du\Big\rangle\,dx - \int_{B_R}\Big\langle
A_\e(x,Du), Du_\eps-Du\Big\rangle\,dx\cr\cr &=&\int_{B_R}\Big\langle
A_\eps(x,Du_\eps), Du_\eps-Du\Big\rangle\,dx -\int_{B_R}\Big\langle
A(x,Du), Du_\eps-Du\Big\rangle\,dx \cr\cr && \quad
-\int_{B_R}\Big\langle A_\eps(x,Du)-A(x,Du),
Du_\eps-Du\Big\rangle\,dx
\end{eqnarray}

\medskip

\noindent Using $\varphi=u$ and $\varphi=u_\e$ as test functions in
\eqref{variationalinequality2}  and \eqref{variationalinequality}
respectively we have

$$\int_{B_R}\Big\langle A_\eps(x,Du_\eps),
Du_\eps-Du\Big\rangle\,dx\le 0 \qquad {\rm and} \qquad
-\int_{B_R}\Big\langle A(x,Du), Du_\eps-Du\Big\rangle\,dx \le 0,$$
therefore from the inequality \eqref{convforte} we deduce

\begin{eqnarray}\label{convforte2}
&&\int_{B_R}(\mu^2+|Du|^2+|Du_\eps|^2)^{\frac{p-2}{2}}|Du_\eps-Du|^2\,dx
\le -\int_{B_R}\Big\langle A_\eps(x,Du)-A(x,Du),
Du_\eps-Du\Big\rangle\,dx\cr\cr &\le &\left(\int_{B_R}|
A(x,Du)-A_\eps(x,Du)|^{\frac{p}{p-1}}\,dx\right)^{\frac{p-1}{p}}\left(\int_{B_R}|Du-Du_\eps|^p\,dx\right)^{\frac{1}{p}}.
\end{eqnarray}

\medskip
Since $p\ge 2$, by well known means, from previous inequality, we
deduce
$$ \int_{B_R}|Du-Du_\eps|^p\,dx\le \int_{B_R}| A(x,Du)-A_\eps(x,Du)|^{\frac{p}{p-1}}\,dx$$
Taking the limit as $\eps\to 0$ in previous inequality, by virtue of
\eqref{convdf}, we deduce that $u_\eps$ converges strongly to $u$ in
$W^{1,p}(B_R)$ and therefore a.e. in $B_R$ for a not relabeled
sequence.

\noindent The strong convergence of $u_\eps$ to  $u$ in
$W^{1,p}(B_R)$ implies also that $u_\eps$ converges strongly to $u$
in $L^{p^*}(B_R)$ and allows us to pass to the limit in \eqref{final
estimate_eps }. So that, by virtue of the Fatou's Lemma and
\eqref{convk} ,  we get
\begin{align*}
\int_{B_{\frac{r}{4}}}\left|DV_p(Du_\e(x))\right|^2 dx\le
\frac{c(\Arrowvert \psi\Arrowvert_{L^{\infty}}^2+\Arrowvert
u_\e\Arrowvert_{L^{p^*}\left(B_{r}\right)}^2)}{r^\frac{p+2}{2}} &
\int_{B_{r}}\left[1+\left|D^2\psi\right|^\frac{p+2}{2}+
\left|D\psi\right|^\frac{p+2}{2}+\kappa^{p+2}+\left|Du\right|^p
\right] dx.
\end{align*}
i.e. the conclusion.
\end{proof}
\noindent Acknowledgements: The third author has been partially supported by the Gruppo Nazionale per l'Analisi Matematica, la Probabilit\`{a} e le loro Applicazioni (GNAMPA) of the Istituto Nazionale di Alta Matematica (INdAM) and by Universit\`{a} degli Studi di Napoli Parthenope through the projects \textquotedblleft sostegno alla Ricerca individuale\textquotedblright(triennio 2015 - 2017) and \textquotedblleft Sostenibilit\`{a}, esternalit\`{a} e uso efficiente delle risorse ambientali\textquotedblright(triennio 2017-2019).

\noindent {\bf M. Caselli}\\
\noindent ETH Z\"{u}rich, Department of Mathematics\\
\noindent {\em E-mail address}: mcaselli@student.ethz.ch

\bigskip
\bigskip

\noindent {\bf A. Gentile}\\
Universit\`{a} degli Studi di Napoli ``Federico II'' Dipartimento di
Mat.~e Appl. ``R.~Caccioppoli'', Via Cintia, 80126 Napoli, Italy

\noindent {\em E-mail address}: andrea.gentile@unina.it

\bigskip
\bigskip

\noindent {\bf R. Giova}\\
\noindent Universit\`{a} degli Studi di Napoli ``Parthenope'' \\
Palazzo Pacanowsky - Via Generale Parisi, 13 \\
 80132 Napoli, Italy

\noindent {\em E-mail address}: raffaella.giova@uniparthenope.it


\begin{thebibliography}{99}

\bibitem{BCGOP} A. L. Baison, A. Clop, R. Giova, J. Orobitg, A. Passarelli di
Napoli, {\em Fractional differentiability for solutions of nonlinear
elliptic equations.}  Potential Anal. 46 (2017), no. 3, 403--430.
\bibitem{BDM}V. B\"{o}gelein, F. Duzaar, G. Mingione, {\em Degenerate problems with irregular obstacles. }J. Reine Angew. Math. 650 (2011), 107--160. 

\bibitem{BFM}M. Bildhauer, M. Fuchs, G. Mingione, {\em A priori gradient bounds and local $C^{1,\alpha}$-estimates for (double) obstacle problems under non-standard growth conditions.} Z. Anal. Anwendungen 20 (2001), no. 4, 959--985.

\bibitem{BLZ} S.-S. Byun, S. Liang, S.-Z. Zheng, {\em  Nonlinear gradient estimates for
double phase elliptic problems with irregular double obstacles}.
Proceedings of the AMS, to appear.

\bibitem{BLOP} S.-S. Byun, K.-A. Lee, J. Oh, J. Park, {\em  Regularity results of the
thin obstacle problem for the p(x)-Laplacian}. J. Functional
Analysis 276 (2019), 496-519.


\bibitem{CKP}{ M. Carozza, J. Kristensen, A. Passarelli di Napoli}, {\em Higher differentiability of minimizers of convex variational integrals.}
Annales Inst. H. Poincar\'e (C) Non Linear Analysis ,  {\bf 28}
(2011), no. 3, 395--411.

\bibitem{CEP} {M. Caselli, M. Eleuteri, A. Passarelli di Napoli}, {\em Regularity
 results for a class of obstacle problems with $p,q$- growth
 conditions} arXiv:1907.08527

\bibitem{Choe} H. J. Choe, {\em A Regualrity Theory for a General Class of Quasilinear
Elliptic Partial Differential Equations and Obstacle Problems.}
Arch. Rational Mech. Anal. 114 (1991),  383-394.

\bibitem{CL} H. J. Choe, J. L. Lewis, {\em On the obstacle problem for quasilinear
elliptic equations of p-Laplace type.} SIAM J. Math. Anal. 22
(1991), no. 3, 623-638.

\bibitem{CGP} A. Clop, R. Giova, A. Passarelli di Napoli, {\em Besov
regularity for solutions of p-harmonic equations}.  Adv. Nonlinear
Anal. 8 (2019), no. 1, 762-778.

\bibitem{EHL} M. Eleuteri, P. Harjulehto, T. Lukkari, {\em Global regularity and
stability of solutions to obstacle problems with nonstandard
growth}. Rev. Mat. Complut. 26  (2013),(1), 147-181.

\bibitem{EleMarMas} {M. Eleuteri, P. Marcellini, E. Mascolo:}\emph{ Lipschitz estimates for systems with ellipticity conditions at infinity,}
Ann. Mat. Pura e Appl. (4),  195 (2016)   1575-1603.

\bibitem{EP}{M. Eleuteri, A. Passarelli di Napoli}, {\em Higher differentiability for solutions to a class of obstacle problems.}
 Calc. Var. Partial Differential Equations 57 (2018), no. 5,  115,
 29 pp.

 \bibitem{EP2}{M. Eleuteri, A. Passarelli di Napoli}, {\em Regularity
 results for a class of non-differentiable obstacle problems} DOI: 10.1016/j.na.2019.01.024

\bibitem{Fuchs}M. Fuchs, {\em H\"{o}lder continuity of the gradient for degenerate variational inequalities}. Nonlinear Anal. 15 (1990), no. 1, 85--100.

\bibitem{FM} {M. Fuchs, G. Mingione}, {\em Full $C^{1,\alpha}$-regularity for free and
constrained local minimizers of elliptic variational integrals with
nearly linear growth}. Manuscripta Math. 102 (2000), 227-250.

\bibitem{Gavioli}{C. Gavioli} {\em Higher differentiability of solutions to a class of obstacle problems under non-standard growth
conditions}. Forum Mathematicum
DOI:https://doi.org/10.1515/forum-2019-0148.

\bibitem{Gentile} {A. Gentile}, {\em Regularity for minimizers of non-autonomous non-quadratic functionals in the case $1<p<2$: an a priori estimate}, Rend. Acc. Sc. fis. mat. Napoli, Vol LXXXV (2018) 185--200.

\bibitem{GiannettiPassa}F. Giannetti, A. Passarelli di Napoli, {\em Higher differentiability of minimizers of variational integrals with variable exponents.} Math. Z. 280 (2015), no. 3-4, 873--892.

\bibitem{Giova1}  R. Giova, {\em Higher differentiability for n-harmonic systems with Sobolev coefficients. }
  J. Differential Equations 259 (2015), no. 11,  5667--5687.


\bibitem{Giova2}  R. Giova, {\em  Regularity results for non-autonomous functionals with $L\log {L}$ -growth and Orlicz Sobolev coefficients.}
 NoDEA Nonlinear Differential Equations Appl.  23  (2016),  no. 6, Art. 64, 18 pp.
 
 \bibitem{GP}
 R. Giova, A. Passarelli di Napoli,
 {\em Regularity results for a priori bounded minimizers of non
 	autonomous functionals with discontinuous coefficients.} Adv. Calc.
 Var. 12 (2019), no. 1, 85-110.
 
\bibitem{Giusti}{E.~Giusti},
{\em Direct methods in the calculus of variations}. World
Scientific, 2003.


\bibitem{Hajlasz} P. Hajlasz. {\em
	Sobolev Spaces on an Arbitrary Metric Space}, Potential Anal.  {\bf
	5} (1996), 403--415.

\bibitem{KristensenMingione} J. Kristensen and G. Mingione, {\em Boundary Regularity in Variational Problems}, Arch. Rational Mech. Anal. 198 (2010) 369--455.
\bibitem{KuusiMingione} T. Kuusi and G. Mingione, {\em Universal potential estimates}, Journal of Functional Analysis 262 (2012) 4205-4269.timates, {\em Journal of Func-\indent tional Analysis}, {\bf 262}, 4205-4269.

\bibitem{Lindqvist}P. Lindqvist, {\em Regularity for the gradient of the solution to a nonlinear obstacle problem with degenerate ellipticity}. Nonlinear Anal. 12 (1988), no. 11, 1245--1255.

\bibitem{MaZ} L. Ma, Z. Zhang {\em Higher differentiability for
solutions of nonhomogeneous elliptic obstacle problems} J. Math.
Anal. Appl. 479 (2019), no. 1, 789-816.

\bibitem{MZ} J. Mu, W. P. Ziemer, {\em Smooth regularity of solutions of double
obstacle problem involving degenerate elliptic equations}. Comm. PDE
16 (1991), nos.4-5, 821-843.


\bibitem{APdN1}  A.~Passarelli di Napoli,
{\em Higher differentiability of minimizers of
 variational integrals with Sobolev coefficients}.
Adv.~Cal.~Var.  { 7}  (2014), no. 1,   59--89.


\bibitem{APdN2} A.~Passarelli di Napoli,
{\em Higher differentiability of solutions of
 elliptic systems with Sobolev coefficients: the case $p=n=2$.} Pot.~Anal.  {\ 41} (2014), no. 3,
 715--735.
\end{thebibliography}
\end{document}